\documentclass{amsart}

\usepackage{amsmath, amssymb, amsfonts, amsthm}
\usepackage{xspace}
\usepackage{mathtools}
\usepackage[all, 2cell]{xy}
\usepackage{tikz}
\usepackage{tikz-cd}
\usetikzlibrary{matrix}
\usetikzlibrary{arrows}
\usetikzlibrary{positioning}
\usetikzlibrary{decorations.pathmorphing}

\theoremstyle{plain}
\newtheorem{td}{Theorem-Definition}

\newtheorem{main}{Theorem}[td]

\newtheorem{theorem}{Theorem}[section]
\newtheorem{lemma}[theorem]{Lemma}
\newtheorem{corollary}[theorem]{Corollary}
\newtheorem{proposition}[theorem]{Proposition}

\newcommand{\bR}{\mathbb{R}}
\newcommand{\bQ}{\mathbb{Q}}

\newcommand{\bZ}{\mathbb{Z}}
\newcommand{\bA}{\mathbb{A}}
\newcommand{\bI}{\mathbb{I}}
\newcommand{\bF}{\mathbb{F}}

\newcommand{\fC}{\mathfrak{c}}
\newcommand{\fF}{\mathfrak{f}}
\newcommand{\fP}{\mathfrak{p}}
\newcommand{\fA}{\mathfrak{a}}
\newcommand{\cO}{\mathcal{O}}
\newcommand{\cU}{\mathcal{U}}

\newcommand{\bmid}{\:\mid\:}

\newcommand{\Set}[1]{{\left\{\,#1\,\right\}}}
\newcommand{\GRH}{\textnormal{GRH}\xspace}
\newcommand{\BPI}{\mathrm{\bf BPI}\xspace}
\newcommand{\BPV}{\mathrm{\bf BPV}\xspace}

\newcommand{\CL}{\mathrm{C\ell}}
\newcommand{\Cl}{\mathrm{c\ell}}

\newcommand{\DIV}{{\,\mathop{divides}\,}}

\newcommand{\GL}{\mathrm{GL}}

\begin{document}

\title{Grunwald-Wang Theorem, an effective version}
\date{}
\author{Song Wang}
\address{Song Wang, Academy of Mathematics
 and Systematics Sciences, the Morningside Center of Mathematics,
 the Key Laboratory of Hua, Chinese Academy of Science.}
\thanks{The author was supported in part by national 973 project 2013CB834202,
also National Natural Science Foundation of China (NNSFC 11321101), and also by the
One Hundred Talent's Program from Chinese Academy of Science.}
\email{\texttt{songw1973@amss.ac.cn}}

\begin{abstract}
The main purpose of this note is to establish an effective version
of the Grunwald--Wang Theorem, which asserts that given a family of
local characters $\chi^{v}$ of $K_{v}^{*}$ of exponent $m$ where
$v \in S$ for a finite set $S$ of primes of $K$, there exists a
global character $\chi$ of the idele class group
$C_{K}$ of exponent $m$ (unless some
special case occurs, when it is $2 m$) whose component at $v$ is
$\chi^{v}$. The effectiveness problem for this theorem is to bound
the norm $N (\chi)$ of the conductor of $\chi$ in terms of $K$, $m$,
$S$ and $N (\chi^{v})$.

The Kummer case (when $K$ contains $\mu_{m}$) is easy since it is almost
an application of the Chinese Remainder Theorem. In this note,
we solve this problem completely in general case,
and give three versions of bound, one is with \GRH, and the other
two are unconditional bounds.

These effective results have some interesting applications in
concrete situations. To give a simple example, if we fix $p$ and
$l$, one gets a good least upper bound for $N$ such that $p$ is not
an $l$--th power mod $N$. One also gets the least upper bound for
$N$ such that $l^{r} \mid \phi (N)$ and $p$ is not an $l$--th power
mod $N$.

Some part of this note is adopted (with some revision)
from my unpublished thesis (2001)
(\cite{Wang2001}).
\end{abstract}
\maketitle

\tableofcontents

\newpage

\section{Introduction} \label{S:1}

Let $K$ be a number field, $K_{v}$ the local completion
of $K$ at $v$, with multiplicative group $K_{v}^{\times}$.
$\bA_{K}$ the ring of adeles of $K$, $\bI_{K} = \bA^{\times}_{K}$
the group of ideles, and $C_{K} = \bI_{K} / K^{\times}$ the group
of idele classes.
\medskip

In 1933, W.\@ Grunwald (\cite{Grunwald33}) stated a striking theorem
asserting that, given a finite set
$S$ of places of $K$, and a family of characters $\chi^{v}$ ($v \in S$)
of $K_{v}^{\times}$ of orders $m_{v}$
(which are also called \emph{local characters}), there exists
a continuous character $\chi$ of finite order of $C_{K}$
(which is also called a \emph{global character})
whose local component
$\chi_{v}$ at each $v \in S$, which are obtained by composing
$\chi$ with the natural local embedding $K_{v}^{\times} \to C_{K}$,
is exactly $\chi^{v}$.
Furthermore, the order of $\chi$ can be made as
the least common multiple $m$ of the $m_{v}$.

\medskip

However, Grunwald's original statement and proof had a flaw,
which occurred when he discussed the special case.
The gap was filled by Sh.\@ Wang (\cite{ShWang48}, \cite{ShWang50}),
who also gave a precise criterion for the special case
and showed that the
order of $\chi$ can be taken to be $m$ under an additional condition
(see Section ~\ref{S:2} and also Chapter 10, \cite{A-T68}).
So this theorem is appropriately called
the \emph{Grunwald-Wang theorem}.

\medskip

Given the local datum $\Set{\chi^{v} \bmid v \in S}$,
there are infinitely many global characters $\chi$
with local components $\chi^{v}$.
They can be highly ramified in general.
Indeed, some additional ramification has to be
allowed as seen for example even when all the $\chi^{v}$
are unramified and $K$ has class number $1$. (If we remove the
order requirement of $\chi$, then there exists a $\chi$ with
given local behavior, and more over, $\chi$ is unramified
outside all finite places outside $S$.)
However, the natural question which arises is whether we can
control the ramification of $\chi$
in terms of $K$, $S$, $m$ and the norms $N(\chi^{v})$
of the conductors of $\chi_{v}$ for $v$ in $S$.

\medskip

Unfortunately, the original proof by Grunwald and
Wang was not effective,
and did not answer this question.
This effectiveness problem was first encountered in 1995 in the work
of J.\@ Hoffstein and D.\@ Ramakrishnan on Siegel zeros (\cite{H-R95}),
where they needed it in a particular case when $K \supset \mu_{m}$ {\bfseries and} $m$
a prime (the Kummer case). As roughly (not precisely) each character
of $K_{v}$ (or $\bI_{K} = \bA_{K}^{\times}$) of exponent $m$ corresponds to a Kummer
extension $K_{v} (\sqrt[m]{y_{v}}) / K_{v}$ or $K (\sqrt[m]{y}) / K$ by the Kummer
theory, the problem can be converted directly into the situation
when the Chinese Remainder Theorem applies (see \cite{Wang2001}).
Although such method is extendable to also the general Kummer case
($K$ contains $\mu_{m}$)
and also the case when $K$ is general and $m$ is an odd prime, it fails
in the general case. So the question in general is not easy.

\medskip

In the thesis \cite{Wang2001}, we solved this problem completely.
we are going to explain our main results. Without loss of generality,
may assume that $m = l^{r}$ a prime power since
it is easy for us to reduce to this case since
any local/global character of finite
order can be written uniquely as a product of characters of prime power degree.

\medskip

Now we introduce some notations. For each place $v$ of $K$,
let $\fP_{v}$ be the symbol of the formal prime of $K$ corresponding to $v$,
which can be used also as the prime idea corresponding to $v$ of $\cO_{K}$
or $\cO_{K_{v}}$,
the ring of integers of $K$ or $K_{v}$ respectively  (depending on the context around).
Let $\chi_{v}$ ($\chi$ resp.\ )
be a continuous character of $K_{v}^{\times}$ ($C_{K}$ resp.\ ). Define
the \emph{arithmetic conductor} or \emph{norm} as the following:

\[
N (\chi_{v}) = \begin{cases}
1 &\text{($\chi_{v}$ is unramified or $v$ is real or complex)} \\
q_{v}^{n} &\text{($n$ is the smallest integer such that
   ${(1 + \fP_{v}^{n})}^{\times} \subset {\rm Ker} (\chi_{v})$)}
\end{cases}
\]
where $q_{v}$ is $1$ if $v$ is archimedean and
the size of residue field of $K_{v}$ when $v$ is finite.
\[
N (\chi) = \prod_{v} N (\chi_{v})
\]

\medskip

Moreover, $N_{S} = \prod_{v \in S} q_{v}$, and $n_{K} = [K : \bQ]$.
Also denote $S_{\infty}$ be the set of infinite places of $K$.
Moreover, for each finite set $T$ of places of $K$ containing
all infinite ones, denote $\CL_{K, T}$ the $T$-class group of $K$
and $\Cl_{K, T}$ its size.

\medskip

By convention, for two quantities $f$ and $g$ depending on various parameters,
$f << g$ or $f = O (g)$ means
$f < C g$ for some absolute positive constant, and
$f <<_{v, v', \ldots} g$ or $f = O_{v, v', \ldots} (g)$
means $f < C g$ for some constant $C$ only depending on $v, v', \ldots$
which are parameters shown up in the subscript.

\medskip


\begin{td} \label{TM:A}
Let $K$ be a number field, $m = l^{r}$ a prime power,
$S$ a finite set of places of $K$,
$\chi^{v}$  a local character of $K_{v}^{\times}$
of exponent $m$ is given for each $v \in S$.
Assume $K(\zeta_{l^{r}})/K$ is cyclic, or the special case occurs
with some specified condition holds \textnormal{(See Section \ref{S:2})}.
There is a global character $\chi$ of $C_{K}$ of exponent $m$
with its local parts $\chi_{v} = \chi^{v}$
for all $v \in S$.

\medskip

Moreover, denote $\BPI$ the smallest number of $N (\chi)$
for such $\chi$.
\end{td}

\medskip

Of course, $\BPI < \infty$ as when $K, S, m = l^{r}$
is given, the number of the choice of $\{ \chi^{v}, v \in S \}$
is finite.

\medskip

In \cite{Wang2001} (unpublished), we show that
\[
\BPI \leq A
{\prod_{v \in S} N(\chi^{v})}^{B}
\]
with $A = {(A_{0}
N_{S})}^{C_{1} |S|}$ where $A_{0}$, $C_{1}$ are constants only depending
on $K$ and $m$, not on $S$, and moreover for the special case
when $K \supset \mu_{m}$
\textbf{OR} $m$ is a prime
We also proved for a better bound
$A = A_{0} N_{S}^{C_{2}}$,
where $A_{0}$, $B$, $C_{1}$ and
$C_{2}$ are independent of $S$. The later
bound improves the result of \cite{H-R95}. The bound $\BPI$ in
this paper is improved compared with \cite{Wang2001}.

\medskip

\begin{main}  \label{TM:1}
All notations are as in Theorem ~\ref{TM:A}. We have
\[
\log \BPI << l^{r} (|S \cup S_{\infty}| + D) \log (d_{K} N_{S} l^{n_{K} r})
\]
where $D = \gamma_{l} (\CL_{K, S \cup S_{\infty}})$, the minimal size of generating
set of the Sylow $l$-subgroup of $\CL_{K, S \cup S_{\infty}}$, which is also
the $\bF_{l}$-rank of $\CL_{K, S \cup S_{\infty}} [l]$,
$l$-torsion part of the class group $\CL_{K, S \cup S_{\infty}}$.
\end{main}

\medskip

\begin{main}  \label{TM:2}
All notations as in Theorem ~\ref{TM:A} and Theorem ~\ref{TM:1}. We have
\[
\BPI <<_{\epsilon, K, l, r} N_{S}^{E_{1} (1/2 + \epsilon)}
\prod_{v \in S} N (\chi^{v})
\]
where
\[
E_{1} = [K (\zeta_{l^{r}}) : K] (|S \cup S_{\infty}| + D - \delta') + \delta
\]
and $D = \gamma_{l} (\CL_{K, S \cup S_{\infty}})$ and
\[
\delta = \begin{cases}
0 &\textnormal{($K (\zeta_{l^{r}}) = K$)} \\
1 &\textnormal{($K (\zeta_{l^{r}}) \ne K$)}
\end{cases}
\]
\[
\delta' = \begin{cases}
1 &\textnormal{($\zeta_{l} \notin K$.)} \\
0 & \textnormal{(otherwise.)}
\end{cases}
\]
\end{main}

\medskip

For example, when $K = \bQ$, $m = l^{r}$ and $S = \{p\}$,
we have $\CL_{K, S \cup S_{\infty}} = 1$ and hence $D = 0$. Moreover,
$\delta' = 0$ if $l = 2$ and $1$ if $l$ is odd; $\delta = 0$ if $l^{r} = 2$,
$1$ if $l^{r} = 4$ or odd prime power and $2$ if $l = 2$ and $r \geq 3$.
\[
E_{1} = \begin{cases}
l^{r - 1} (l - 1) + 1 & \text{($l$ is odd.)} \\
2 & \text{($l^{r} = 2$.)} \\
2^{r} + 1 & \text{($l = 2$, $r \geq 2$.)}
\end{cases}
\]
(Note that if $2 \notin S$ and $r \geq 3$, then we can replace $\delta = 1$ by $0$.)

\medskip

The following theorem is a result with \GRH (General Riemann Hypothesis)
(which is an analogue of Riemann Hypothesis to the case of general $L$-functions).

\medskip

\begin{main}  \label{TM:3}
All notations as in Theorem ~\ref{TM:A} and Theorem ~\ref{TM:1}. Also assume
GRH holds. Then there is an effectively computable $C > 0$ only depending
on $l, r$ such that
\[
\BPI \leq {(C \log (d_{K} N_{S} l^{n_{K}}))}^{2 (|S \cup S_{\infty}| + D - \delta' + \delta)}
\prod_{v \in S} N (\chi^{v})
\]
where $D = \gamma_{l} (\CL_{K, S \cup S_{\infty}})$, $\delta$
and $\delta'$ are defined in Theorem ~\ref{TM:1} and Theorem ~\ref{TM:2}.
\end{main}

\medskip

\emph{Remark}: $\delta = 1$ can be replaced by
$0$ in the following two cases: (1) $\zeta_{l^{r}} \notin K$, $K (\zeta_{l^{r}}) / K$ is cyclic,
and $S_{l} \not\subset S$, (2) $l = 2$ and $K (\zeta_{2^{r}}) / K$ is not cyclic
and $S_{0} \not\subset S$.

\medskip

\emph{Remark}: In fact, the quantity $|S \cup S_{\infty}| + D - \delta'$
is exactly the $\bF_{l}$-rank of the Selmer group ${\rm Sel}_{l}^{S} (K)$.

\medskip

These effective results have some interesting applications in
concrete situations. To give a simple example, if we fix $p$ and
$l$, one gets a good least upper bound for $N$ such that $p$ is not
an $l$--th power mod $N$. One also gets the least upper bound for
$N$ such that $l^{r} \mid \phi (N)$ and $p$ is not an $l$--th power
mod $N$.

\medskip


Let's explain the idea of these theorems right now. There are two main ingredients.
First, we are trying to reduce this effective problem to a variant. Roughly
speaking, Let $P = \prod_{v \in S} K_{v}^{\times} \subset \bI_{K}$ and
$P_{0} = {\rm Ker} \prod_{v} \chi^{v} \subset P$. To find the global $\chi$
with given local behavior, it suffices to find a suitable standard
open subgroup $V$ of $\bI_{K}$ such that
$$K^{\times} \bI_{K}^{m} V \cap P \subset P_{0} \qquad(*)$$
This is exactly how the Grunwald-Wang was proved in \cite{A-T68}.
To find the least bound for $N (\chi)$, it suffices to
choose $V = \prod_{v} V_{v}$ carefully to ensure (*)
 (Of course for $v \in S$ we choose $V_{v}$ lying in the kernel of $\chi^{v}$).
Let $\gamma i^{m} u \in P$ where $i \in \bI_{K}$, $\gamma \in K^{\times}$
and $u = (u_{v}) \in V$. Then $m | v (\gamma)$ for $v \notin S$ finite.
Then $\gamma = \gamma' c^{m}$ for some $\gamma'$ lies in
certain finite generate abelian subgroup $A$ of $K^{\times}$.
Choose a finite set $T$ of finite places of $K$, disjoint from $S$,
such that $A \cap \bigcap_{v \in T} (1 + \fP_{v}) \subset A^{m}$.
Now let $V_{v} = 1 + \fP_{v}$ when $v \in T$, and $V_{v} = \cO_{v}$
outside $T \cup S$, and hence such $V$ meets our need.
We are able to do that since $A$ is finite generated.
(In fact, $A [l]$, the $l$-torsion part of such $A$ is isomorphic
to the Selmer group
\[
{\rm Sel}_{l}^{S} (K) = \{ [x] \in K^{S} / {K^{S}}^{l},
x \in {K_{v}^{\times}}^{l} \,\,\forall v \in S\}
\]
)

\medskip

Now, pick any $\alpha \in A - A^{m}$, we can choose $v$ outside $S$
such that $v$ is inert at $K (\zeta_{m}, \sqrt[m]{\alpha}) / K (\zeta_{m})$
so that $\alpha \notin {K_{v}^{\times}}^{m}$. A series $v$ can ensure
that $A \cap \bigcap_{v} (1 + \fP_{v}) \subset A \cap {K (\zeta_{m})^{\times}}^{m}$.
Then add at most two more $v$ to make the intersection to $A^{m}$. Over $v$s
$K (\zeta_{m}, \sqrt[m]{A}) / K$ does not collapse, i.e.,
$\forall x \in K (\zeta_{m}, \sqrt[m]{A})$, if $x \in K_{v}$ for all $v \in T$
then $x \in K$ (See section ~\ref{S:2}).

\medskip

Now the second ingredient comes. To find $v$ and bound the norm as above,
we are always in the following situation: Let $E/F$ be a Galois extension,
and we want to fine the least prime ideal $\fP$ of $F$ outside a given finite
set $S$ not split in $E$. Some $S$-version of Chebotarev Density Theorems
(more precisely, $S$-version of Multiplicity One Theorems
in $GL (1)$, see Section ~\ref{S:3} and ~\ref{S:4})
will answer this question. Our three main theorem rely
on exactly three $S$-versions of multiplicity one theorems
(also see Theorem ~\ref{T:311}, Theorem ~\ref{T:401}).

\medskip

\begin{theorem} \label{T:101}
\textnormal{\textbf{($S$-version of the Multiplicity One for $GL (1)$, Theorem ~\ref{T:401})}}

\medskip

Let $K$ be a number field and $\chi$ a nontrivial global character of $C_{K}$ of finite order.
Then there is a place $v$ of $K$ such that

\medskip

\textnormal{(1)} $\fP_{v} \notin S$.

\medskip

\textnormal{(2)} $\chi_{v} \ne 1$ and is not ramified.

\medskip

\textnormal{(3 A)} $\log N (\fP_{v}) << \log A (\chi, S)$

\medskip

\textnormal{(3 B)} $N (\fP_{v}) <<_{\epsilon, K} N (\chi)^{1/2 + \epsilon} N_{S}^{\epsilon}$
for every $\epsilon > 0$.

\medskip

\textnormal{(3 C)} With \GRH, $N (\fP_{v}) << {(\log A (\chi, S))}^{2}$.

\medskip

where $A (\chi, S) = d_{K} N (\chi) N_{S}$.

\end{theorem}

\medskip

In fact, (A) and (C) are $S$-effective versions of \cite{L-M-O79} and \cite{L-O77},
and the proof runs through some variant of the proof of the above two papers
(For details, see \cite{Wang2013-2}).
(B), uses Landau's method, and runs arguments like \cite{LW2009}
with some modifications. We include the proof in Section ~\ref{S:4}.

\medskip

This paper is organized as the following:
In Section ~\ref{S:2},
after the notation, we introduce the Grunwald-Wang Theorem, and the
special case of Wang. Experts can skip most parts of this section.
In Section ~\ref{S:3} we will formulate and prove the variant and hence our theorems.
As mentioned before, three $S$-version of multiplicity one for $GL (1)$
will be formulated there and will be proved in Section ~\ref{S:4}.

\medskip

Most parts of Section \ref{S:3} and some parts of Section \ref{S:4}
are adopted from Chapter 5 of my thesis \cite{Wang2001} with certain revision.
Theorem \ref{TM:2} is new since it adopts results from other people in later publications
(\cite{Brumley06}, \cite{Wyh2006}, \cite{LW2009}).
Here we express again our great appreciation for
my advisor D. ~Ramakrishnan for his suggestion of this topic,
and for the continued help during my graduate study and continuing years.

\bigskip

\section{Preliminaries} \label{S:2}


In this part, we'll recall the Grunwald-Wang Theorem,
and also make some preliminaries.
First we recall some notations.

\medskip

Again, $K$ denotes a number field or a function field in this section,
$S$ a finite set of places of $K$.
For each place $v$, $K_{v}$ denotes the completion of $K$
at $v$ which is a complete discrete valuation field.
When $v$ is non-archimedean, denote $\cO_{K_{v}}$ the valuation
ring of $K_{v}$ and $\cU_{K_{v}} = \cO_{K,_{v}}^{\times}$
the unit group of $K_{v}$.

\medskip

Denote $\bA_{K}$ the ring of Adeles of $K$, $\bI_{K} = \bA_{K}^{\times}$
the group of ideles of $K$.
Denote $C_{K} = \bI_{K} / \times K^{\times}$
the idele class group of $K$.
When $S$ contains $S_{\infty}$, the set of all archimedean places
of $K$ which is empty when $K$ is a function field,
 denote $\bA_{K, S}$ and $\bI_{K, S} = \bA_{K, S}^{\times}$
the ring of $S$-adeles and the group of $S$-ideles respectively.
Recall that $\bA_{K, S} = \prod_{v \in S} K_{v} \times \prod_{v \notin S} \cO_{K_{v}}$
and $\bA_{K} = \bigcup_{S \supset S_{\infty}} \bA_{K, S}$. $K$ and $K^{\times}$
embed into $\bA_{K}$ and $\bI_{K}$ diagonally respectively.
Denote $K^{S} = K^{\times} \cap \bI_{K, S}$ the set of
$S$-units (which are those $a \in K$ such that
$a \in \cU_{K_{v}})$ for $v \notin S$.

\medskip


\subsection{Grunwald-Wang}

\hspace{1ex}

\medskip

First recall the following results which is crucial for the statement and the proof of
the Grunwald-Wang Theorem. Denote $\zeta_{n}$ a primitive root
of unity of the order $n$,
and we choose $\zeta_{n}$ for all $n$ to ensure that $\zeta_{n m}^{m} = \zeta_{n}$.
Denote $\eta_{n} = \zeta_{n} + \zeta_{n}^{-1}$.

\medskip

\begin{proposition} \label{T:201}
\textnormal{\textbf{(\cite{A-T68}, Theorem 10.1)}}

Let $P(m, S)$ be the set of all $x \in K^{\times}$ such that
$x \in {(K_{v}^{\times})}^{m}$
for all $v \notin S$.

\medskip

Then $P(m, S) = {K^{\times}}^{m}$, unless exactly if the following
conditions hold. \textnormal{(If they hold, we say that \emph{the special case
of Wang} occurs.)}

\medskip

\textnormal{(a)} $K$ is a number field.

\medskip

\textnormal{(b)}  $-1, \pm (2 + \eta_{2^{s}})$ are non-squares in $K^{\times}$
where $s$ is the integer such that $\eta_{2^{s}} \in K^{\times}$,
$\eta_{2^{s + 1}} \notin K^{\times}$.

\medskip

\textnormal{(c)} $m = 2^{t} m'$, $2 \nmid m'$, $t \gneqq s$,

\medskip

\textnormal{(d)} $S_{0} \subseteq S$, where $S_{0}$ is the set of places
$v | 2$ such that  $-1, \pm(2 + \eta_{2 ^{s}})$
are nonsquares in $K_{v}$.

\medskip

In the special case, $P(m, S) = {K^{\times}}^{m} \cup a_{0} {K^{\times}}^{m}$,
where $a_{0} = (1 + \zeta_{2 ^{s}})^{m} = \eta_{2^{s+1}}^{m}
= (i \eta_{2^{s+1}})^{m} = [\pm (2 + \eta_{2^{s}})]^{m / 2}$.
\end{proposition}

\medskip

\emph{Remark}: When the special case of Wang occurs,
$F = K (\sqrt{-1}, \sqrt{2 + \eta_{2^s}})$
is a $F_{4} = \bZ / 2 /bZ \times \bZ / 2 \bZ$ extension of $K$,
and for all $v$ such that $K_{v} (\sqrt{-1}, \sqrt{2 + \eta_{2^s}}) / K_{v}$
is also an $F_{4}$-extension, $v \in S_{0}$. In particular, $t \geq 3$.

\medskip

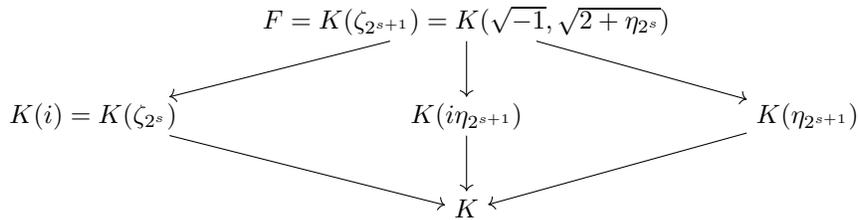
\[
\begin{tikzcd}
& F = K (\zeta_{2^{s + 1}}) = K (\sqrt{-1}, \sqrt{2 + \eta_{2^s}})
\arrow{dl}{} \arrow{d}{} \arrow{dr}{} & \\
K (i) = K (\zeta_{2^{s}}) \arrow{dr}{}
& K (i \eta_{2^{s + 1}}) \arrow{d}{}
& K (\eta_{2^{s + 1}}) \arrow{dl}{} \\
& K &
\end{tikzcd}
\]

\medskip

\emph{Remark}: The simplest special case of Wang is the following: $K = \bQ$, $m = 8$,
$S = \{2\}$ and $a_{0} = 16$. One simple example of the special case of Wang
when $S = \infty$ is the following: $K = \bQ (\sqrt{7})$, $m = 8$ and $a_{0} = 16$.

\medskip

Now we state the first version of the Grunwald-Wang Theorem. Recall that a character
is said to be of exponent $m$ if its order divides $m$.

\medskip

\begin{theorem} \label{T:202}
\textnormal{(\cite{A-T68}, Theorem 10.5-1)}

\medskip

Let $S$ be a finite set of places of $K$,
and $\chi^{v}$, a local character of $K_{v}^{\times}$
of exponent $m$ for each $v \in S$.

\medskip

Then there is a global character $\chi$ of $C_{K}$ of exponent $m$
with its local parts $\chi_{v} = \chi |_{K_{v}^{\times}} = \chi^{v}$
for all $v \in S$ unless the special case of
Wang \textnormal{(Proposition ~\ref{T:201})}
occurs.

\medskip

When the special case of Wang occurs, such $\chi$ can still be chosen of exponent $m$
when the following condition holds

\[
\prod_{v \in S} \chi_{v} (a_{0}) = 1
\]
where $a_{0}$ is defined as in
\textnormal{Proposition ~\ref{T:201}}.
Otherwise, such $\chi$ can be chosen of exponent $2 m$.
\end{theorem}

\medskip

By the class field theory, the local and global version, we have the following variants
which is the corollary.

\medskip

\begin{theorem} \label{T:203}
\textnormal{(\cite{A-T68}, Theorem 10.5-2)}

\medskip

Let $S$ be a finite set of places of $K$,
and $L^{v} / K_{v}$ a cyclic extension
of degree $m_{v}$ for each $v \in S$.

\medskip

Then there is a cyclic extension $L / K$ of degree $m = {\rm LCM} (m_{v})$,
the least multiple of $m_{v}$ for $v \in S$,
such that $L_{w} \cong L^{v}$ for the place $w$ of $L$ above $v$ for $v \in S$
unless the special case of Wang \textnormal{(Proposition ~\ref{T:201})}
occurs.

\medskip

When the special case of Wang occurs, such $L$ can still be chosen of degree $2 m$
over $K$.
\end{theorem}

\medskip

In next section we will formulate a variant.
It is known that this variant implies Theorem ~\ref{T:202} (\cite{A-T68}).
We will prove this variant in an effective way.

\medskip


\subsection{Estimations}

\hspace{1ex}

\medskip

In this part, we list some results of estimations.

\medskip

\begin{proposition}  \label{T:204}
\textnormal{\textbf{(Discriminant-Conductor Formula)} (\cite{L-O77},
 \cite{Odlyzko77})}

\medskip

Let $L / K$ be an abelian extension of local (resp.\ number fields)
 of degree $n$, and $\chi_{0} = 1$,
$\chi_{1}$, \ldots, $\chi_{n - 1}$ be all $n$ local
character (resp.\ idele class characters)
of $K^{\times} / N_{L / K} (L^{\times})$ $C_{K} / N_{L / K} C_{L}$. Then

\[
d_{L} = d_{K}^{n} \prod_{i = 1}^{n-1} N (\chi_{i})
\]

\end{proposition}

\medskip

\qedsymbol

\medskip

\emph{Remark}: By the class field theory, ${\rm Gal} (L / K) \cong
K^{\times} / N_{L / K} L^{\times}$ (resp.\ ${\rm Gal} (L / K) \cong
C_{K} / N_{L / K} C_{L}$) when $L / K$ is an abelian extension
of local fields (resp.\ number fields) of degree $n$,
and hence there are exactly $n = [L : K]$
characters.

\medskip

\begin{lemma}  \label{T:205}

\medskip

Let $L / K$ be an abelian extension of local fields (resp.\ number fields)
of prime degree $l$, and $\chi$ a local character (resp.\ idele class character)
associated to $L / K$. Then

\[
d_{L} = d_{K}^{l}  {N (\chi_{i})}^{l - 1}
\]

\end{lemma}

\medskip

\emph{Proof}: This is a corollary to the previous lemma. \qed

\medskip

Recall that $N_{S}$ denotes the set of norms of (formal) primes
$\fP_{v}$ for all $v \in S$. Recall that
\[
N (\fP_{v}) = \begin{cases}
1 &\text{($v$ is real and complex.)} \\
q_{v} = \# (\cO_{K_{v}} / \fP_{v}) &\text{($v$ is finite.)}
\end{cases}
\]

\medskip

Now we define the \emph{conductor cycle} $\fF (\chi_{v})$
of a local character $\chi_{v}$ of $K_{v}^{\times}$ (of finite order)
as following:
\[
\fF (\chi_{v}) = \\
\begin{cases}
1 &\text{($\chi_{v} = 1$, $v$ is real or complex.)} \\
\fP_{v} &\text{($\chi_{v} = {\rm Sgn}$, $v$ is real. )} \\
1 &\text{($\chi_{v}$ is unramified, $v$ is finite.)} \\
\fP_{v}^{n} &\text{($\chi_{v}$ is ramified, and $n$ is the smallest)}\\
 &\text{( integer such that $\chi_{v} ( 1 + \fP_{v}^{n}) = 1$)}.
\end{cases}
\]
and $N (\chi_{v}) = N (\fF (\chi_{v}))$.

\medskip

A \emph{formal cycle} of $K$ is a formal product $\prod_{v} \fP_{v}^{n_{v}}$
where $n_{v} = 0, 1$ when $v$ is real, $n_{v} = 0$ when $v$ is complex,
$n_{v} \geq 0 \in \bZ$ for other $v$, and $n_{v} = 0$ for almost all $v$.
Define the \emph{norm} of
$\fF = \prod_{v} \fP_{v}^{n_{v}} = \prod_{v} {(N (\fP_{v}))}^{n_{v}}$.

\medskip

For a idele class character $\chi$ of finite order, denote
its \emph{conductor cycle} $\fF (\chi) = \prod_{v} \fF (\chi_{v})$
where $\chi_{v}$ is the restriction of $\chi$ to $K_{v}^{\times}$.
Denote $\fF_{\rm f} (\chi) = \prod_{v \text{finite}} \fF (\chi_{v})$.
Define the \emph{arithmetic conductor} or the \emph{norm} of $\chi$
as $N (\chi) := N (\fF (\chi))$, the norm of its conductor cycle,
and the \emph{analytic conductor} or \emph{level} $A (\chi)$ of $\chi$
as $A (\chi) := d_{K} N (\chi)$.

\medskip

\begin{lemma} \label{T:206}
Let $K_{v}$ be a local field with valuation $v$ and assume
that $p$ is the characteristic of its residue field. Then
if $v (x) > v (m) + \frac{e_{v / p}}{p - 1}$ where $e_{v / p} = v (p)$
the ramification index of $K_{v} / \bQ_{p}$,
then $1 + x \in {\cO_{K_{v}}^{\times}}^{m}$. In addition, if $p \nmid m$
and $v (x) > 0$,
then $1 + x \in {\cO_{K_{v}}^{\times}}^{m}$.
\end{lemma}

\medskip

\emph{Proof}: Want to prove that under this condition
the binomial series ${(1 + x)}^{1 / m} = \sum_{k = 0}^{\infty} \binom{1 / m}{k} x^{k}$
converges. Let $a_{n} = \binom{1 / m}{k} x^{k}$ where
\[
\binom{N}{k}
 = \frac{N (N - 1) \cdots (N - k + 1)}{k!}
\]

\medskip

\begin{align}
v (a_{k}) &= v (\frac{1}{m} (\frac{1}{m} - 1) \cdots \frac{1}{m} + k + 1)
  - v (k!) + k v (x)
\notag \\
&\geq k (v (x) - v (m)) - v (k!) \notag \\
&\geq k (v (x) - v (m) - e_{v / p} (p-1)) \notag \\
&\to 0 \qquad (k \to \infty)\notag
\end{align}

\medskip

Here the second inequality holds since the multiplicity of the factor $p$ in $k!$
is $[k / p] + [k / p^{2}] + \cdots$ which is not greater than
$k / p - 1$.

\medskip

For the case when $p \nmid m$, the assertion follows since
we can prove that $\binom{1 / m}{k}$ a $p$-adic integer.

\medskip

\qed

\medskip

\begin{lemma} \label{T:207}

\medskip

\textnormal{(1)} If $\chi$ is a global character of
$C_{K}$ of order $m$, then
\[
\fF_{\rm f} (\chi) \DIV (m) \prod_{p | m} (p) \prod_{v \in S_{\rm Ram}}
\fP_{v}
\]
where $S_{\rm Ram}$ is the set of finite places $v$
where $\chi_{v} = \chi |_{K_{v}^{\times}}$ is ramified. Moreover,
\[
N (\chi) \leq {\left( m \prod_{p | m} p \right)}^{n_{K}} N_{S_{\rm Ram}}
\]
where $n_{K} = [K : \bQ]$.

\medskip

\textnormal{(2)} In particular,
if $m = l^{r}$ where $r > 0$, then
\[
N (\chi) \leq l^{n_{K} (r + 1)} N_{S_{\rm Ram}}
 \leq l^{n_{K} (r + 2)} N_{S_{\rm Ram} - S_{l}}
\]
where $S_{l}$ is the set of places $v$ of $K$ such that $v | l$.
\end{lemma}

\medskip

\emph{Proof}:

\medskip

(1) Assume $\fF_{\rm f} (\chi) = \prod_{v \in S_{\rm Ram}} \fP_{v}^{n_{v}}$.
From Lemma ~\ref{T:206}, $n_{v} = 1$ if $v \nmid m$.
If $p | m$ and $v | p$, then
$n_{v} \leq 1 + v (m) + e_{v / p} / (p - 1) \leq 1 + v (m) + v_{p} = 1 + v (p m)$.
Thus

\begin{align}
\fF_{\rm f} (\chi) & = \prod_{v \in S_{\rm Ram}} (\chi)
\prod_{p | m} \prod_{v | p} \fP_{v}^{v (p m)} \notag \\
& \DIV \prod_{v \in S_{\rm Ram}} (\chi) (m) \prod_{p | m} (p) \notag
\end{align}

\medskip

Take the norm, we get the last assertion in (1). (2)
is routine to check when we take $m = l^{r}$.

\medskip

\qedsymbol

\medskip

\emph{Remark}: In Lemma ~\ref{T:207} (2), when $l$ is odd,
then we can have
\[
N (\chi) \leq l^{n_{K} r} N_{S_{\rm Ram}}
\]
instead. This follows from Lemma ~\ref{T:206} (2) also by checking carefully.

\medskip

\begin{lemma}  \label{T:208}
\textnormal{(1)} Let $K' = K (\zeta_{l^{r}})$. Then
\[
d_{K'} \leq  {(l^{n_{K} (r + 2)})}^{[K' : K] - 1} d_{K}^{[K' : K]}
\]

\medskip

\textnormal{(2)} Let $\tilde{\chi}$ be a global character
of $C_{K'}$ of order $l$ associated to the Kummer extension
$K' (\sqrt[l]{a}) / K'$ via the global class field theory
where $a \in K^{\times}$. Let $T_{a}$ be the set
of finite places $v$ of $K$ such that $v \nmid l$
and $l \nmid v (a)$. Then
\[
N (\tilde{\chi}) \leq l^{n_{K'} (r + 2)} N_{T_{a}}^{[K' : K]}
= {(l^{n_{K} (r + 2)} N_{T_{a}})}^{[K' : K]}
\]
and
\[
A (\tilde{\chi}) \leq
{(l^{n_{K} 2 (r + 2)} d_{K} N_{T_{a}})}^{[K' : K]}
\]

\medskip

\textnormal{(3)} Let $\chi'$ be a global character
of $C_{K}$ of order $l$ associated to the subextension
of the cyclic extension $K (\zeta_{l^{r}}) / K$.
Then
\[
N (\chi') \leq l^{3 n_{K}}
\]
and
\[
A (\chi') \leq l^{3 n_{K}} d_{K}
\]
\end{lemma}

\medskip

\emph{Proof}: Let $n' = [K' : K]$ where $K' = K (\zeta_{l^{r}})$.

\medskip

(1): Let $\chi_{0} = 1$, $\chi_{1}$, $\ldots$, $\chi_{n' - 1}$ be the $n'$
characters of $C_{K} / N_{L / K} C_{L}$. From Proposition ~\ref{T:204},
the Discriminant-Conductor formula, we have
$d_{K'} \leq D_{K}^{n'} \prod_{i = 1}^{n - 1} N (\chi_{i})$.
Now all assertions follow from Lemma ~\ref{T:206} (2),
since all $\chi_{i}$
are ramified at finite places $v \nmid l$, as by the class field theory,
the Galois extension $K'_{i} / K$ associated to $\chi_{i}$
are all subextensions of the cyclotomic extensions
$K' = K (\zeta_{l^{r}}) / K$.

\medskip

(2): Note that the following fact: For each $a \in K^{\times}$,
the Kummer extension
$K' (\sqrt[l]{a}) / K'$ is ramified at $v \nmid l$ if and only if
$l \nmid v (a)$, hence $\tilde{T}'$, the set of places
of $K'$ not dividing $l$ ramify in this Kummer extension,
is contained in the set of places $w$ of $K'$ above $T_{a}$.
Hence $\tilde{T}' \leq N_{T_{a}}^{[K' : K]}$.

\medskip

Now it is routine to check (2) from Lemma ~\ref{T:207} (1) and (2).

\medskip

(3): Clear from Lemma ~\ref{T:207} (2).

\medskip

\qedsymbol

\bigskip

\section{The Variant and the Proof} \label{S:3}

In this section, we will formulate a variant of the Grunwald-Wang Theorem.
Also, we'll prove this variant and then the main theorems of this
paper.

\medskip

\subsection{The Variant}

\hspace{1ex}


\medskip

Now $K$ again denotes a number field, and $S$
a finite set of places of $K$.
Let $P = \prod_{v \in S} K_{v}^{\times} \hookrightarrow \bI_{K}$
endowed with the induced topology from $\bI_{K}$,
$P_{0}$ a subgroup of $P$ containing $P^{m}$.
When the special case of Wang occurs,
denote $\mathbf{\rm a}_{0} = (a_{0}, \ldots a_{0}) \in P$
where $a_{0} \in K^{\times}$ is defined as in
Proposition ~\ref{T:201}.

\medskip

\begin{theorem} \label{T:301}
\textnormal{(\textbf{\cite{A-T68}, Theorem 10.4})}

\medskip

Let $S$, $P$, $P_{0}$, $\mathbf{a}_{0}$ as above. Then their exists an open subgroup
$V$ of $\bI_{K}$ such that
$K^{\times} \bI_{K}^{n} V \cap P \subset P_{0}$
where $n$ is an integer such that
$K^{\times} \bI_{K}^{n} \cap P \subset P_{0}$.
Moreover, $n$ can be taken as $m$ unless the special case of Wang
occurs. When the special case of Wang holds,
then $n$ can be taken as $m$ when $\mathbf{a}_{0} \in P_{0}$.
Otherwise, $n$ be taken as $2 m$.
\end{theorem}

\medskip

\qedsymbol

\medskip

\emph{Remark}: As proven in \cite{A-T68},
Theorem ~\ref{T:301} implies the Grunwald-Wang
since when given $\{\chi^{v}, v \in S\}$, we can set
$P_{0} = {\rm Ker} \prod_{v} \chi^{v}$,
and $N = P_{0} \bI_{K}^{m} K^{\times} V$, then
$P / P_{0} \cong \bI_{K} / N$ through the natural embedding
$i_{P}: P \hookrightarrow \bI_{K}$. Thus $\chi$
can be chosen as
\[
\chi: \bI_{K} / K^{\times} \to \bI_{K} / N
\overset{i_{P}^{-1}}{\longrightarrow} P / P_{0}
\overset{\prod_{v \in S} \chi^{v}}{\longrightarrow}
\mathbb{C}^{\times}
\]

\medskip

Before we formulate the effective version of the variant, we recall
some definitions. Let $\fC = \prod_{v} \fP_{v}^{n_{v}}$
be a formal cycle of $K$. Define the \emph{standard
open subgroup} $V_{\fC}$ as $\prod_{v} \cU_{v}^{(n_{v})}$
where
\begin{align}
\cU_{v}^{(0)} = K_{v}^{\times} &\text{($v$ is real or complex)}
\notag \\
\cU_{v}^{(1)} = \bR^{+, \times} &\text{($v$ is real)}
\notag \\
\cU_{v}^{(0)} = \cU_{K_{v}} &\text{($v$ is finite)}
\notag \\
\cU_{v}^{(n)} = 1 + \fP_{v}^{n} \subset \cU_{K_{v}}
&\text{($v$ is finite, $n \geq 1$)} \notag
\end{align}

\medskip

\begin{td} \label{TM:B}
Let $K$ be a number field, $m = l^{r}$ a prime power,
$S$ a finite set of places of $K$,
$P = \prod_{v \in S} K_{v}^{\times} \hookrightarrow \bI_{K}$,
$P_{0}$ a subgroup of $P$ containing $P^{m}$.
Assume that $K^{\times} \bI_{K}^{m} \cap P \subset P_{0}$.

\medskip

Then there is a standard open subgroup
$V_{\fC}$ of $\bI_{K}$ such that
$K^{\times} \bI_{K}^{m} V_{\fC} \cap P \subset P_{0}$.

\medskip

Denote $\BPV$ be the smallest number of $N (\fC)$
for such $V_{\fC}$.
\end{td}

\medskip

\qedsymbol

\medskip

The main theorems of this paper (Theorem ~\ref{TM:1}, ~\ref{TM:2},
~\ref{TM:3}) follow from the following
theorems.

\begin{theorem} \label{T:302}
All notations are same as in Theorem-Definition ~\ref{TM:A}
and ~\ref{TM:B}. Then $\BPI \leq \BPV$.
\end{theorem}

\medskip

\begin{theorem}  \label{T:303}
All notations are as in Theorem ~\ref{TM:B}. We have
\[
\log \BPV << l^{r} (|S \cup S_{\infty}| + D) \log (d_{K} N_{S} l^{n_{K} r})
\]
where $D = \gamma_{l} (\CL_{K, S \cup S_{\infty}})$, the minimal size of generating
set of the Sylow $l$-subgroup of $\CL_{K, S \cup S_{\infty}}$.
\end{theorem}

\medskip

\begin{theorem}  \label{T:304}
All notations as in Theorem ~\ref{TM:B}, Theorem ~\ref{T:303}. We have
\[
\BPV <<_{\epsilon, K, l, r} N_{S}^{E_{1} (1/2 + \epsilon)} N (\fC_{0})
\]
where
\[
E_{1} = [K (\zeta_{l^{r}}) : K] (|S \cup S_{\infty}| + D - \delta') + \delta
\]
where $D = \gamma_{l} (\CL_{K, S \cup S_{\infty}})$ and
\[
\delta = \begin{cases}
0 &\textnormal{($K (\zeta_{l^{r}}) = K$)} \\
1 &\textnormal{($K (\zeta_{l^{r}}) \ne K$)}
\end{cases}
\]
\[
\delta' = \begin{cases}
1 &\textnormal{($\zeta_{l} \notin K$.)} \\
0 & \textnormal{(otherwise.)}
\end{cases}
\]
and $\fC_{0}$ is the smallest cycle such that $V_{\fC} \cap P \subset P_{0}$.
\end{theorem}

\medskip

\emph{Remark}: When $P / P_{0}$ is cyclic, and $P_{0} =
\prod_{v \in S} \chi^{v}$, then $\fC_{0} = \prod_{v \in S} \fF (\chi^{v})$.

\medskip

\emph{Remark}: $\delta = 1$ can be replaced by
$0$ in the following two cases: (1) $\zeta_{l^{r}} \notin K$, $K (\zeta_{l^{r}}) / K$ is cyclic,
and $S_{l} \not\subset S$, (2) $l = 2$ and $K (\zeta_{2^{r}}) / K$ is not cyclic
and $S_{0} \not\subset S$. This also applies to Theorem ~\ref{T:305}.

\medskip

\begin{theorem}  \label{T:305}
All notations as in Theorem ~\ref{TM:B}, Theorem ~\ref{T:303}
and Theorem ~\ref{T:304}. Also assume
GRH holds. Then there is an effectively computable $C > 0$ only depending on
$l, r$ such that
\[
\BPV \leq {(C \log (l^{n_{K}} d_{K} N_{S}))}^{2 (|S \cup S_{\infty}| + D - \delta' + \delta)}
N (\fC_{0})
\]
where $D = \gamma_{l} (\CL_{K, S \cup S_{\infty}})$, $\delta$,
$\delta'$ and $\fC_{0}$ are defined in Theorem ~\ref{TM:B} and Theorem ~\ref{T:304}.
\end{theorem}

\medskip

Theorem ~\ref{T:302} follows from the remark after
Theorem-Definition ~\ref{TM:B}.
The rest of this paper will be devoted to the proof of these theorems.

\medskip

\subsection{Main Arguments and Proofs}

\medskip

\hspace{2ex}

\medskip

In this part, we will prove Theorem ~\ref{TM:B}, ~\ref{T:303},
~\ref{T:304}, ~\ref{T:305} and finally the main theorem.
Again recall the notation
and our situation. $K$ denotes a number field,
$m = l^{r}$ a prime power, $S$ a finite
set of finite places. Let $K' = K (\zeta_{l^{r}})$.
Recall that $P = \prod_{v \in S} K_{v}^{\times}$,
and $P \supset P_{0} \supset P^{m}$. To make the discussion easier,
we may assume that $K^{\times} \bI_{K}^{m} \cap P \subset P_{0}$
(which means that the special case of Wang occurs and the special condition
to guarantee the order $m$ also holds, see Theorem ~\ref{T:301}).
The goal is to find the least bound $\BPV$ of $N (\fC)$ for $V_{\fC}$
such that $K^{\times} \bI_{K}^{m} V \cap P \subset P_{0}$
(cf.\ Theorem ~\ref{T:301}).

\medskip

\emph{Definition}: Let $E / F$ be an abelian extension
of number fields. Let $T$ be a finite set of places
of $F$. Denote $B (E / F, T, l)$ be the least bound $C$
such that, for each subextension $E'$ of $E / F$
cyclic of degree $l$, there is a prime $\fP_{v}$ of $F$
satisfying: (a) $v \notin T$. (b) $\fP_{v}$ is unramified and
does not split in
$E'$. (c) $N (\fP_{v}) \leq C$.

\medskip

\emph{Remark}: By the Chebotarev Density Theorem, $B (E / F, T, l) < \infty$.

\medskip

\emph{Definition}: Define $P^{*} (n, S)$ as the set of
$x \in K^{\times}$ such that $m | v (x)$ for each $v \notin S$.

\medskip

Moreover, recall that for a finite abelian group $H$, $\gamma_{l} (H)$
denotes the minimal generators of the Sylow-$l$ group

\medskip

\begin{proposition}  \label{T:306}
Let $E / F$ be an abelian $l$-extension of number fields,
$T$ a finite set of places of $K$. Let
$e = \gamma_{l} ({\rm Gal} (E / F))$. Then There
are $e$ places $w_{1}, w_{2}, \ldots, w_{e}$ of $F$
such that

\medskip

\textnormal{(1)} $w_{i} \notin T$ for all $1 \leq i \leq e$.

\medskip

\textnormal{(2)} For any $a \in E - F$,
their is a $w_{i}$ which does not split
in $F (a) / F$.

\medskip

\textnormal{(3)} $N (\fP_{w_{i}}) \leq B (E / F, T, l)$.

\end{proposition}

\medskip

\emph{Proof}: We proceed in steps.

\medskip

(Step 1): Assume first that $E / F$ is an elementary $l$-extension,
i.e., ${\rm Gal} (E / F) \cong \bF_{l}^{e}$.
Change (2) to the following (2)':

Induct in $e$.
$e = 1$ is OK. Assume that the proposition holds for $< e$, and now want
to prove the $e$ case.

\medskip

Pick a cyclic subextension $E_{1} / F$, which is cyclic of degree $l$,
and pick $w_{1} \notin T$ of degree $1$, where $E_{1} / F$ is inert at $w_{1}$,
and moreover $N (\fP_{w_{1}}) \leq B (E_{1} / F, T, l) \leq B (E / F, T, l)$.

\medskip

Let $E_{2}$ be the fixed field of $D_{w_{1}} (E_{1} / F)$, the decomposition group
of $E / F$ at $w_{1}$. Then $[E_{2} : F] = l^{e - 1}$, and $w_{1}$ splits
in $E_{2} / F$ completely. By induction, there are $w_{2}, w_{3}, \dots, w_{e}$
satisfying (1), (2), (3) with $E$ replaced by $E_{2}$.

\medskip

Now, (1) and (3) hold for original $E$ and $w_{1}, w_{2}, \ldots, w_{e}$.
For (2), for each $a \in E - F$, if $a \in E_{2}$,
then (2) for $E'$ and $w_{2}, \dots, w_{e}$
implies (2). If $a \notin E_{2}$, then $D_{w_{1}} (E / F)$ does not fix
$a$, and hence $w_{1}$ does not split in $F (a) / F$ then.

\medskip

(Step-2): Assume that $E / F$ is an $l$-extension. Let $E'$ be the composition
of elementary $l$-subextensions of $E / F$. Thus ${\rm Gal} (E' / F) \cong \bF_{l}^{e}$.
Let $w_{1}, w_{2}, \dots, w_{e}$ be chosen to satisfying (1), (2), (3) for $E'$.
Then (1), (3) for $E$ holds. For (2) for $E$, let $a \in E - F$, want to find
$w_{i}$ such that $F (a) / F$ is not split at $w_{i}$. In fact, the maximal elementary
$l$-subextension of
$F (a) / F$ is not split at $w_{i}$ for some $i$,
and hence $F (a) / F$ is not split at this $w_{i}$.

\medskip

\qedsymbol

\medskip

\begin{proposition} \label{T:307}
Let $K' = K (\zeta_{l^{r}})$. Let $e = |S \cup S_{\infty}| + D - \delta'$.
Assume that $K' / K$ is cyclic.

\medskip

Then there exist places $v_{1}, v_{2}, \ldots, v_{e + \delta}$
of $K$ satisfying the following:

\medskip

\textnormal{(1)} $v_{1}, v_{2}, \ldots, v_{e + \delta} \notin S$.

\medskip

\textnormal{(2)} For any $y \in P^{*} (l^{r}, S)$ such that
$y \in {K_{v_{i}}^{\times}}^{l^{r}} \quad (1 \leq i \leq e + \delta)$,
we have $y \in {K^{\times}}^{l^{r}}$.


\medskip

\textnormal{(3)} $N (\fP_{v_{i}}) \leq $ $
B (K' (\sqrt[l^{r}]{P^{*} (l^{r}, s)}) / K, \tilde{S}, l)$
for $1 \leq i \leq e$ and \\
$N (\fP_{v_{e + 1}}) $ $\leq B (K' / K, S, l)$
if $\delta = 1$.

\medskip

Here $\tilde{S}$ the set of places of $K'$ above the places in $S$,
$D = \gamma_{l} (\CL_{K, S \cup S_{\infty}})$ and
\[
\delta = \begin{cases}
0 &\textnormal{($K (\zeta_{l^{r}}) = K$)} \\
1 &\textnormal{($K (\zeta_{l^{r}}) \ne K$)}
\end{cases}
\]
\[
\delta' = \begin{cases}
1 &\textnormal{($\zeta_{l} \notin K$.)} \\
0 & \textnormal{(otherwise.)}
\end{cases}
\]
\end{proposition}

\medskip

\emph{Remark}: This is the effective version of Proposition ~\ref{T:201}.

\medskip

\emph{Proof}: Recall $m = l^{r}$.

\medskip

(Step-1): Let $\tilde{K} = K' (\sqrt[l^{r}]{P^{*} (l^{r}, S)})$.
Then $\gamma_{l} ({\rm Gal} (\tilde{K} / K')) \leq e$.

\medskip

Let $T$ be a set of places containing $S \cap S_{\infty}$
such that $\CL_{K, T}$ is not divisible by $l$. We can choose
$T$ such that $|T| = |S \cup S_{\infty}| + D$ where
$D = \gamma_{l} (\CL_{K, S \cup S_{\infty}})$. In particular,
if $l \nmid \CL_{K, S \cup S_{\infty}}$ then $D = 0$.
In fact, let
$Q$ be the Sylow $l$-group of $\CL_{K, S \cup S_{\infty}}$
and $c_{1}, \ldots, c_{D}$ be a set of generators of $Q$.
By the Chebotarev density, there exist places $w_{1}, \ldots, w_{D}$
of $K$ such that the ideal class $\fP_{w_{j}}$ in
$\CL_{K, S \cup S_{\infty}}$
is $c_{j}$ for each $1 \leq j \leq D$.

\medskip

Let $T = S \cup S_{\infty} \cup \{w_{1}, \ldots w_{D}\}$.
Then $P^{*} (m, S) \subset {K^{\times}}^{m} K^{T}$. In fact,
for each $y \in P^{*} (m, S)$,
let $\fA_{y} = \prod_{v \notin S \cap S_{\infty}} \fP_{v}^{v (y) / l^{r}}$.
Then the ideal $\fA_{y}^{l^{r}} = (y) J$ for some $J \in J_{K}^{S}$ where
$J_{K}^{S}$ is the ideal group generated
$\fP_{v}$ for $v \in S$.
Hence the class of $\fA_{y}$ in $\CL_{K, S \cup  S_{\infty}}$
has power $l^{r}$. Hence by the choice of $T$,
the class of $\fA_{y}$ in $\CL_{K, T}$ is trivial. In particular,
there is $z \in K^{\times}$ such that $\fA_{y} = (z) \fA'$
where $\fA'$ is some ideal in $J_{K}^{T}$.
Thus $y z^{-l^{r}} \in K^{T}$.
Hence $\tilde{K} \subset \tilde{K}' = K' (\sqrt[l^{r}]{K^{T}})$.
Thus
\[
\gamma_{l} (\tilde{K} / K') \leq
\gamma_{l} (\tilde{K}' / K') \leq {\rm rank}_{\bF_{l}} K^{T} / {(K^{T})}^{l}
\]
Now (Step-1) is done since ${\rm rank} (K^{T}) = |S \cup S_{\infty}| - 1 + D$
(and hence $\tilde{K} / K'$ is a finite abelian extension, and
\[
{\rm rank}_{\bF_{l}} (K^{T} / {(K^{T})}^{l}) = e = |S \cup S_{\infty}| + D - \delta'
\]

\medskip

(Step-2): Assume that $K' = K$, i.e., $\zeta_{l^{r}} \in K$.
Then by Proposition ~\ref{T:306}, there exist places
$w_{1}, \ldots, w_{e}$ of $K$ of outside $S$,
such that, for each $b \in \tilde{K} - K$, $K (b) / K$
is not split at $w_{i}$ for some $i$, and moreover
$N (\fP_{w_{i}}) \leq B (\tilde{K} / K, S, l)$ (by (Step-1)).

\medskip

Now consider any $y \in P^{*} (m, S)$ such that $y \in {K^{\times}}^{m}$.
Set $b \in \sqrt[l^{r}]{y}$, and by the discussion
above, $K (b) / K$ is not split for some $w_{i}$,
and hence $b \notin K$. Hence $y \notin {K^{\times}}^{m}$.
Hence $v_{1} = w_{1}, \ldots, v_{e} = w_{e}$ are what we want.

\medskip

(Step-3): $K' \ne K$ and $\delta = 1$. By Proposition ~\ref{T:306},
there exist places
$w_{1}, \ldots, w_{e}$ of $K'$ of outside $\tilde{S}$,
such that, for each $b \in \tilde{K} - K'$, $K' (b) / K'$
is not split at $w_{i}$ for some $i$, and moreover
$N (\fP_{w_{i}}) \leq B (\tilde{K} / K', \tilde{S}, l)$ (by (Step-1)).
Moreover, there exists a place
$v'_{1}$ of $K$
outside $S$, such that $K' / K$
is inert at $v'_{1}$, and moreover
$N (\fP_{v'_{i}}) \leq B (K' / K, S, l)$ (Proposition ~\ref{T:306},
noting that fact
${\gamma}_{l} (K' / K) = 1$).

\medskip

Now let $v_{1}, \ldots, v_{e}$ be the places of $K$ under
$w_{1}, \ldots, w_{e}$ and $v_{j + 1} = v'_{j}$
($j \leq 1$). Then (1)
and (3) are clear.

\medskip

For (2), Let $y \in P^{*} (m, S)$ and assume that
for all $i = 1, \ldots e + \delta$, $y \in {K_{v_{i}}^{\times}}^{m} \subset
    {{K'}_{w_{i}}^{\times}}^{m}$.
Hence by the choice of $w_{1}, \ldots, w_{e}$,
we have $y \in {{K'}^{\times}}^{m}$.

\medskip

let $K_{1}$
be the maximal $l$-subextension of $K' / K$. First claim
$y \in {K_{1}^{\times}}^{m}$. In fact,
$y^{[K' : K_{1}]} = N_{K' / K_{1}} (y) \in {K_{1}^{\times}}^{m}$, and
hence $y \in {K_{1}^{\times}}^{m}$ since $l \nmid [K' : K_{1}]$
and $([K' : K_{1}], m) = 1$. As $[K_{1} : K]$ is a cyclic $l$-extension,
then it has only one elementary $l$-subextension $K_{2}$.
By the choice of $v_{e + 1} = v'_{1}$, $K_{2} / K$ is not split at $v'_{1}$
and $K_{1} / K$ is inert at $v'_{1}$. Since $y \in {K_{v'_{1}}^{\times}}^{m}$,
there is a $z \in K'$ such that $z^{m} = y$, and $z \in K_{v'_{1}}$,
and hence $K (z) = K$ as $K (z) / K$ is inert at $v'_{1}$.
Hence $y \in {K^{\times}}^{m}$.

\medskip




\qedsymbol

\medskip

\begin{lemma} \label{T:308}
Let $m = 2^{r}$ and assume that $K' = K (\zeta_{m}) / K$ is not cyclic.
Let $K^{+} = K (\eta_{2^{m}})$. Let $s$ be the integer such that
$\eta_{2^{s}} \in K$ but $\eta_{2^{s + 1}} \notin K$.

\medskip

\textnormal{(1)} Let $y \in K (i)$ or $K (i \eta_{2^{s + 1}})$
such that $y^{m} \in K^{\times}$,
then there is an $j \in \bZ$ such that $y {(1 + \zeta_{2^{s}})}^{j} \in K$.

\medskip

\textnormal{(2)} Denote $a_{0} = {(1 + \zeta_{2^{s}})}^{m}$. Assume that $v$
is a place such that $\eta_{2^{s + 1}} \notin K_{v}$,
and $Y \in {K_{v}^{\times}}^{m} \cap K^{\times} \cap {K'^{\times}}^{m}$.
Then $Y \in {K^{\times}}^{m} \cup a_{0} {K^{\times}}^{m}$.

\medskip

\textnormal{(3)} Let $v | 2$ be such that $-1$ and $\pm (2 + \eta_{2^{s}})$
is not square in $K_{v}^{\times}$. Assume $Y \in {K_{v}^{\times}}^{m} \cap K^{\times}
\cap {K'^{\times}}^{m}$
then $Y \in {K^{\times}}^{m}$.
\end{lemma}

\medskip

\emph{Proof}: In this case $s < r$. Moreover, $[K' : K^{+}] = 2$.
Let $\sigma$ be the nontrivial $K^{+}$-automorphism of $K'$.
Then $\sigma (\zeta_{2^{r}}) = \zeta_{2^{r}}^{-1}$.
Moreover $K (\zeta_{2^{s + 1}}) / K$ is a $4$-extension
with three intermediate subextension $K (i)$, $K (\eta_{2^{s + 1}})$
and $K (i \eta_{2^{s + 1}})$.

\medskip

For more
analysis on this case, see Section ~\ref{S:3}.

\medskip

(1) Put $y_{r} = y {(1 + \zeta_{2^{s}})}^{r}$. Then $y_{r} / \sigma (y_{r}) =
(y / \sigma (y)) \zeta_{2^{s}}^{r}$.
Want to prove that $y^{2^{s}} = \sigma (y^{2^{s}})$
so that, can pick $r$ such that $\sigma (y_{r}) = y_{r}$ and hence
$y_{r} \in K$ as $K = {K (i)}^{\sigma} = {K (i \eta_{2^{s + 1}})}^{\sigma}$.

\medskip

Note that $T = y / \sigma (y)$ is a root of unity in $K (i)$ or $K (i \eta_{2^{s + 1}})$.
Hence $T = \zeta_{s}^{j}$ for some $j$ in the first case or $T = \pm 1$
in the later case. As $s \geq 2$, in all case we have $T^{s} = 1$. Hence such $j$ is
what we want.

\medskip

(2) By the choice of $v$, $K^{\times} / K$ is inert at $v$.
Then either $K (i)$ or $K (i \eta_{2^{s + 1}})$ is split
at $v$. Let $K_{2}$ be such field, and $w$ a place above $v$.
By comparing the residue degree,  $K' / K_{2}$
is inert at $w$ as they have the same residue degree as $K^{+} / K$.
Thus $Y \in {K_{v}^{\times}}^{m} = {K_{2, w}^{\times}}^{m}$.
Hence there is some $z$ lies in $K' \cap K_{v} = K_{2, w}$ such
that $z^{m} = y$. As $K' / K_{2}$ is inert at $w$,
$z \in K_{2}$. Thus (2) follows from (1) since ${(1 + \zeta_{2^{s}})}^{m} = a_{0}$.

\medskip

(3) There is a $z \in K' \cap K_{v}$ such that $z^{m} = y$.
By the choice of $v$ and the analysis on $K' / K$ at the beginning,
$K' / K$ does not collapse at $v$. Hence $z \in K$ and $Y = z^{m}
\in {K^{\times}}^{m}$.

\qedsymbol

\medskip

\begin{proposition} \label{T:309}
Let $K' = K (\zeta_{2^{r}})$ and $e = |S \cup S_{\infty}| + D$.
Assume that $K' / K$ is not cyclic. Denote $s$, $K^{+}$ and
$a_{0}$ are as in Lemma
~\ref{T:308}. Let $S_{0}$ be the set of places $v | 2$ where
$-1$ and $\pm (2 + \eta_{2^{s}})$ is not a square in $K_{v}^{\times}$.

\medskip

Then there exist places $v_{1}, v_{2}, \ldots, v_{e + 1}$
of $K$ satisfying the following:

\medskip

\textnormal{(1)} $v_{1}, v_{2}, \ldots, v_{e + 1} \notin S$.
If $S_{0} \not\subset S$ then $v_{e + 1} \in S_{0} - S$.

\medskip

\textnormal{(2)} For any $y \in P^{*} (2^{r}, S)$ such that
$y \in {K_{v_{i}}^{\times}}^{2^{r}} \quad (1 \leq i \leq e + 1)$,
we have $y \in {K^{\times}}^{2^{r}}$ if $S_{0} \not\subset S$
and $y \in {K^{\times}}^{2^{r}} \cup a_{0} {K^{\times}}^{2^{r}}$
if $S_{0} \subset S$.

\medskip

\textnormal{(3)} $N (\fP_{v_{i}}) \leq $ $
B (K' (\sqrt[l^{r}]{P^{*} (2^{r}, s)}) / K, \tilde{S}, 2)$
for $1 \leq i \leq \delta$ and \\
$N (\fP_{v_{e + 1}}) $ $\leq B (K' / K, S, 2)$
if $S_{0} \subset S$.

\medskip

Here $\tilde{S}$ the set of places of $K'$ above the places in $S$,
$D = \gamma_{2} (\CL_{K, S \cup S_{\infty}})$.
\end{proposition}

\medskip

\emph{Proof}: Recall $m = 2^{r}$.
Note that (Step-1) of the proof of Proposition
~\ref{T:307} still applies since we don't use the assumption
that $K (\zeta_{l^{r}}) / K$ is cyclic. Moreover,
(Step-3) goes through smoothly except for the verifications
of (2).

\medskip

(Step-3'): Let $v_{1}, \ldots, v_{e}$ be the places of $K$ under
$w_{1}, \ldots, w_{e}$ as in (Step-3) the proof of Proposition
~\ref{T:306} and $v_{j + \delta} = v'_{j}$
($j \leq \delta$) if $S_{0} \subset S$ and any place
in $S_{0} - S$ if $S_{0} \not\subset S$.
Then again (1) and (3) are clear.

\medskip

For (2), Let $y \in P^{*} (m, S)$ and assume that
for all $i = 1, \ldots e + 1$, $y \in {K_{v_{i}}^{\times}}^{m}
    \subset {{K'}_{w_{i}}^{\times}}^{m}$.
Hence by the choice of $w_{1}, \ldots, w_{e}$,
we have $y \in {{K'_{w_{i}}}^{\times}}^{m}$
and hence $y \in {{K'}^{\times}}^{m}$.

\medskip

Now (2) follows from Lemma ~\ref{T:308}.

\qedsymbol

\medskip

\begin{theorem} \label{T:310}
All notations are as in Theorem-Definition ~\ref{TM:B}, Theorem ~\ref{T:301},
Theorem ~\ref{T:302}, Proposition ~\ref{T:307} and Proposition ~\ref{T:309}.
Then
\[
\BPV \leq l^{n_{K} (r + 2)} N (\fC_{0})
B (K' (\sqrt[l]{P^{*} (l^{r}, S)}) / K', \tilde{S}, l)
 ^{|S \cup S_{\infty}| + D - \delta'}
B (K' / K, S, l)^{\delta}
\]
where $\fC_{0}$ is the smallest cycle such that $V_{\fC} \cap P \subset P_{0}$.
\end{theorem}

\medskip

\emph{Proof}: Again denote $e = |S \cup S_{\infty}| + D - \delta'$.
If we are in the special case of Wang, let $\mathbf{a}_{0} = (a_{0}) \in P$
be the idele with component $a_{0}$ at all $v \in S$ and $1$ elsewhere.

\medskip

Let $v_{1}$, $\ldots$, $v_{e + \delta}$ be the set of places chosen
in Proposition ~\ref{T:307} (resp.\ Proposition ~\ref{T:309})
when $K' = K (\zeta_{l^{r}})$ is cyclic (resp.\ $K' / K$ is not cyclic).

\medskip

Now consider three cycles as following: $\fC_{0}$ the conductor
cycle of $P_{0}$ as in Theorem ~\ref{T:305}, Theorem ~\ref{T:306}
and Theorem ~\ref{T:307}. $\fC_{1} = {(l)}^{r + 1} \prod_{v | l} \fP_{v}$
and $\fC_{2} = \prod_{i = 1}^{e + \delta} (\fP_{v_{i}})$
Let $\fC = \fC_{0} \fC_{1} \fC_{2}$, and we have $N (\fC)$ is bounded
by the right hand side of the formula in this theorem (Proposition ~\ref{T:307},
Proposition ~\ref{T:309}).

\medskip

Want to prove $K^{\times} \bI_{K}^{m} V_{\fC} \cap P \subset P_{0}$.
Let $\gamma \in K^{\times}$, $B = (B_{v}) \in \bI_{K}$ and $u = (u_{v}) \in V_{\fC}$
such that $Y = \gamma B^{m} u \in P$. Want $Y \in P_{0}$.

\medskip

In fact, $\gamma B_{v}^{m} u_{v} = 1$ for $v \notin S$, $\gamma \in P^{*} (m, S)$.
Moreover, for each $v = v_{i}$, if $v \nmid l$ (resp.\ if $v_{i} | l$),
then $U_{v} = 1 + \fP_{v} \subset {K_{v}^{\times}}^{m}$ (resp.\
$U_{v} = 1 + \fP_{v}^{(r + 1) v (l) + 1} \subset {K_{v}}^{\times}$
(Lemma ~\ref{T:205}, ~\ref{T:207})).
Hence $\gamma \in {K_{v_{i}}^{\times}}^{m}$ for $i = 1, \ldots e + \delta$.

\medskip

Thus by Proposition ~\ref{T:307} and Proposition ~\ref{T:309},
we have $\gamma \in {K^{\times}}^{m}$, or
$\gamma \in {K^{\times}}^{m} \cup a_{0} {K^{\times}}^{m}$
when we are in the special case of Wang and $S_{0} \subset S$.

\medskip

When $\gamma \in {K^{\times}}^{m}$, then
for $v \in S$, $Y_{v} = \gamma B_{v}^{m} u_{v} \in u_{v} {K_{v}^{\times}}^{m}$,
and hence $Y \in (P \cap V_{\fC_{0}}) P^{m} \subset P_{0}$.

\medskip

When we are in the special case of Wang
and $\gamma \in {K^{\times}}^{m} \cup a_{0} {K^{\times}}^{m}$,
we have for each $v$,
$Y_{v} \gamma B_{v}^{m} u_{v} \in u_{v}
 ({K_{v}^{\times}}^{m} \cup a_{0} {K_{v}^{\times}}^{m})$,
and hence $Y \in (P \cap V_{\fC_{0}}) (P^{m}
    \cup \mathbf{a}_{0} P^{m}) \subset P_{0}$
as $\mathbf{a}_{0} \in P_{0}$.
Done.

\qedsymbol

\bigskip

\subsection{Estimations again}


In this part, we'll plug in all estimations we made in Section ~\ref{S:1}
to finish the proof.

\medskip

We have the following three effective versions of the Chebotarev density and multiplicity one,
which will be proved in next section.

\medskip

For each global character $\chi$ of $C_{K}$, define $A (\chi, S) :=
A (\chi) \prod_{v \in S, \text{$\chi_{v}$ is unramified}} \leq A (\chi) N_{S}$.

\medskip

\begin{theorem} \label{T:311}
There is a prime $\fP_{v}$ satisfying the following:

\medskip

\textnormal{(1)} $\fP_{v} \notin S$.

\medskip

\textnormal{(2)} $\chi_{v} \ne 1$ is unramified.

\medskip

\textnormal{(3)} $N (\fP_{v})$ has the following versions of the least bound:

\medskip

\textnormal{(A)} $\log (N (\fP_{v})) << \log A (\chi, S)$.

\medskip

\textnormal{(B)} $N (\fP_{v}) <<_{\epsilon, K} N (\chi)^{1 /2 + \epsilon} N_{S}^{\epsilon}$.

\medskip

\textnormal{(C)} Assuming \GRH, $N (\fP_{v}) << {(\log A (\chi, S))}^{2}$.
\end{theorem}

\medskip

For the proof, see Section ~\ref{S:4}.

\medskip

\begin{corollary} \label{T:312}
Let $L / K$ be a cyclic extension of degree $m$.
There is a prime $\fP_{v}$ satisfying the following:

\medskip

\textnormal{(1)} $\fP_{v} \notin S$.

\medskip

\textnormal{(2)} $\fP_{v}$ is ramified and does not split in $L / K$.

\medskip

\textnormal{(3)} Let $\chi$ be the global character of $C_{K}$ corresponding to
$L / K$ via the class field theory.
$N (\fP_{v})$ has the exactly three versions of
the least bound as the same in \textnormal{(3)}
in Theorem \ref{T:311}.
\end{corollary}

\medskip

\qedsymbol

\medskip

\begin{lemma} \label{T:313}

\medskip

\textnormal{(1)} Let $K' = K (\zeta_{l^{r}})$ and $\tilde{K} = K' (\sqrt[l^r]{P^{*} (l^{r}, S)})$.
Let $\tilde{S}$ be as in Proposition ~\ref{T:307} and Proposition ~\ref{T:309},
and $\tilde{\chi}'$ a global character associated to a cyclic subextension of $\tilde{K} / K'$.
of order $l$.
Then $A (\tilde{\chi}', \tilde{S}) \leq {(l^{2 n_{K} (r + 2)} d_{K} N_{S})}^{[K' : K]}$.

\medskip

\textnormal{(2)} Let $\chi'$ be a global character associated to
a subextension of $K' / K$ of order $l$.
Then $A (\chi', S) \leq {(l^{3 n_{K}})} d_{K} N_{S}$.
\end{lemma}

\medskip

\emph{Proof}: (1) Any subextension of $K' / K$ is a Kummer extension. So applying
Lemma ~\ref{T:208} (2) by noting that $T_{a} \subset S$. (2) follows from Lemma
~\ref{T:208} (3).

\medskip

\qedsymbol

\medskip

\emph{Proof of Theorem ~\ref{T:303}, Theorem ~\ref{T:304} and Theorem ~\ref{T:305} using
Theorem ~\ref{T:310} and Theorem ~\ref{T:311}}:

\medskip

\emph{Proof}: We proceed in three situation:

\medskip

(1) Theorem ~\ref{T:303}:
\begin{align}
\log &\BPV \leq \log N (\fC_{0})
 + e \log B (\tilde{K} / K', \tilde{S}, l)
 + \delta \log B (K' / K, S, l) \notag \\
& \leq \log (N_{S}) + (n_{K} (r + 2)) \log (l)
+ e [K' : K] \log {l^{2 n_{K} (r + 2)} d_{K} N_{S}} \notag \\
&+ \delta \log {l^{3 n_{K}} d_{K} N_{S}} \notag \\
&<< (e [K' : K] + \delta) \log (l^{r n_{K}} d_{K} N_{S}) \notag \\
&<< l^{r} (|S \cup S_{\infty}| + D) \log (l^{r n_{K}} d_{K} N_{S}) \notag
\end{align}

\medskip

(2) Theorem ~\ref{T:304} and Theorem ~\ref{T:305}: Can be proved similarly
by using Theorem ~\ref{T:310} and Theorem ~\ref{T:311}.

\qedsymbol

\emph{Proof of Theorem ~\ref{TM:1}, Theorem ~\ref{TM:2}, Theorem ~\ref{TM:3}}:

\medskip

\emph{Proof}: Use Theorem ~\ref{T:302}, Theorem ~\ref{T:303}, Theorem ~\ref{T:304},
Theorem ~\ref{T:305}.

\medskip

\qedsymbol

\bigskip

\section{The $S$-effective Version of
  the Strong Multiplicity One -- $GL (1)$ Case} \label{S:4}

In this section, we'll prove Theorem ~\ref{T:311}
and Corollary ~\ref{T:312}.


First we summarize the rough idea of Theorem ~\ref{T:311} (A) and (C). The main idea
of them is the same as in \cite{L-O77} and \cite{L-M-O79}.
Theorem ~\ref{T:311} (A) (without \GRH) and (C) (with \GRH) follows from a variant of Theorem 1.2
and Corollary 1.3 (plus its proof) in \cite{L-M-O79}. The only
difference is that we replace the quantities $d_{L}$, $A (\chi)$ by $d'_{L} = d_{L} N_{S}^{[L : K]}$
and $A (\chi, S) = A (\chi) N_{S}$. More details can be also found in \cite{Wang2013-2}.

\medskip

In this section we prove (B).
which uses Landau's idea
(see \cite{Landau27}, the proof follows almost from \cite{LW2009}).

\medskip

Now we restate Theorem ~\ref{T:311} (A) and (B) as the following theorems.

\medskip

\begin{theorem} \label{T:401}
\textnormal{\textbf{($S$-version of the Multiplicity One for $GL (1)$)}}

\medskip

Let $K$ be a number field and $\chi$ a nontrivial global character of $C_{K}$ of finite order.
Then there is a place $v$ of $K$ such that

\medskip

\textnormal{(1)} $\fP_{v} \notin S$.

\medskip

\textnormal{(2)} $\chi_{v} \ne 1$ and is not ramified.

\medskip

\textnormal{(3 A)} $\log N (\fP_{v}) << \log A (\chi, S)$

\medskip

\textnormal{(3 B)} $N (\fP_{v}) <<_{\epsilon, K} N (\chi)^{1/2 + \epsilon} N_{S}^{\epsilon}$
for every $\epsilon > 0$.

\medskip

\textnormal{(3 C)} With \GRH, $N (\fP_{v}) << {(\log A (\chi, S))}^{2}$.

\medskip

where $A (\chi, S) = d_{K} N (\chi) N_{S}$.

\end{theorem}

\bigskip

\subsection{Proof \textnormal{(B)}}:

\medskip \hspace{2ex} \medskip

In this part, we prove Theorem ~\ref{T:401} (which is Theorem ~\ref{T:311} (B)).
We'll quote a lot of standard results of number theory (for references of textbook style
see \cite{Lang-GTM110}, \cite{Lang-GTM211}, \cite{Neukirch91}, \cite{RV-GTM186},
\cite{Weil74} etc, for related topics,
see \cite{Brumley06}, \cite{H-R95}, \cite{L-M-O79}, \cite{L-O77},
\cite{Wyh2006}, \cite{LW2009}, \cite{RW2003}, \cite{Wang2001}, \cite{Stk74} etc).
To make things easier, may assume that for all $v \in S$, if $v$ is a finite place,
and $\chi_{v}$ is unramified. In this case $A (\chi, S) = A (\chi) N_{S}$.

\medskip


\emph{Proof of Theorem ~\ref{T:401}}:

\medskip

Set
\[
S (X, \chi, S) = \sum_{n = 1}^{\infty} a_{\chi, S} (n) \omega (\frac{n}{X})
\]
where $a_{\chi, S} (n)$ is the $n$-th coefficient of
the Hecke $L$-function $L^{S} (s, \chi) = $ \\
$\prod_{v \notin S} {(1 - \chi_{v} (\pi_{v}) N (\fP_{v})^{-s})}^{-1}$ \\
$= \sum_{n = 1}^{\infty} a_{\chi, S} (n) n^{-s}$,
and the weight function $\omega (x)$ defined as a smooth function which
may be specified as \cite{Wyh2006}:
\[
\omega (X) = \begin{cases}
0 &\text{($x \leq 0$ or $x \geq 3$)} \\
e^{-1/x} &\text{($0 < x \leq 1$)} \\
e^{-1/(3 - x)} &\text{($2 \leq x < 3$)} \\
\leq 1 &\text{(all $x$)}
\end{cases}
\]

\medskip

Consider the Mellin transform
\[
W (s) = \int_{0}^{\infty} \omega (x) x^{s - 1} d\,x
\]
which is a analytic function of $s$. Fix $\sigma < 0$ and let $s = \sigma + i t$
then
\[
W (s) <<_{A, \sigma} \frac{1}{{(1 + |t|)}^{A}}
\]
for all $A > 0$ by repeated partial integration.
By Mellin inversion,
\[
\omega (x) = \frac{1}{2 \pi i} \int_{(2)} W (s) x^{-s} d\, s
\]
where the integration is made along the vertical line ${\rm Re} s = 2$.

\medskip

Then we have
\begin{align}
S (X, \chi, S) &= \frac{1}{2 \pi i} \sum_{n = 1}^{\infty} a_{\chi, S} (n) \int_{(2)}
  W (s) {\left( \frac{n}{X} \right)}^{-s} d\, s \notag \\
&= \frac{1}{2 \pi i} \int_{(2)} X^{s} W (s) L^{S} (s, \chi) \notag
\end{align}
where the interchange of the summation and the integral is guaranteed by the absolute
convergence along the real line $\sigma = 2$. By the standard arguments,
plus the fact that all incomplete $L$-functions of nontrivial characters are entire
of order $1$,
we may shift the integral line to get
\[
S (X, \chi, S) = \frac{1}{2 \pi i} \int_{-H} x^{s} W (s) L^{S} (s, \chi) d\,s
\]
where $H > 0$ is to be specified later.

\medskip

Let $L_{S} (s, \chi) = \prod_{v \in S} L (s, \chi_{v})$,
and $L_{\infty} (s, \chi) = L (s, \chi_{\infty})$ the
gamma factor of $\chi$. Then we have the following
functional equation
\[
L^{S} (s, \chi) = W (\chi) A (\chi)^{(1/2 - s)} G_{0} (s) G_{1} (s)
  L^{S} (1 - s, \chi^{-1})
\]
where $W (\chi)$ is the root number of $\chi$ which has absolute value $1$,
\[
G_{0} (s) = \frac{L_{\infty} (1 - s, \chi^{-1})}{L_{\infty} (s, \chi)}
\]
and
\[
G_{1} (s) = \frac{L_{S} (1 - s, \chi^{-1})}{L_{S} (s, \chi)}
\]
We need to estimate $G_{0} (s)$ and $G_{1} (s)$ along the vertical line
$\sigma = -H$, avoiding to the pole of them.

\medskip

\begin{lemma} \label{T:402}

\medskip

\textnormal{(1)} If $H > 0$ is not an integer,
then
\[
G_{0} (-H + iT) <<_{H, n_{K}} {( 1 + |t|)}^{n_{K} (1/2 + H)}
\]

\medskip

\textnormal{(2)}
$G_{1} (-H + iT) <<_{H, n_{K}, \epsilon'} N_{S}^{H + \epsilon'}$.
\end{lemma}

\medskip

We quote the following results on gamma function for which the proof
can be found in a lot of analysis textbooks.

\medskip

\begin{lemma}  \label{T:403}
\textnormal{\textbf{(Stirling formula)}}

\[
|\Gamma (\sigma + it)| = \sqrt{2 \pi} e^{- \frac{\pi}{2} t} {|t|}^{\sigma - 1/2}
    \left( 1 + O_{\sigma, \delta, A} \left( \frac{1}{1 + |t|} \right) \right)
\]
for all $|t| > A > 0 $ and $s = \sigma + it$ away from any poles by at least distance $\delta$.
\end{lemma}

\medskip

\qedsymbol

\medskip

\emph{Proof of Lemma ~\ref{T:402}}:

\medskip

(1) It is well known that there are $a (\chi)$ and $b (\chi)$ be the integer such that
$L_{\infty} (s, \chi) = C \pi^{-s n_{K}/2}
    {(\Gamma (s / 2))}^{a (\chi)} {(\Gamma ((s + 1) / 2))}^{b (\chi)}$
and $a (\chi) + b (\chi) = n_{K}$. Note that $a (\chi) = a (\chi^{-1})$.
Now assume $H > 0$, $|T| > 1$, then
\begin{align}
G_{0} (-H + iT) &= \frac{{L_{\infty} (1 + H - iT, chi^{-1})}}{L_{s} (-H + iT, \chi)} \notag \\
&= \pi^{n_{K} (-1/2 - H + iT)} \frac{\Gamma ((1 + H - iT) / 2)}{\Gamma ((-H + iT) / 2)}
\frac{\Gamma ((H - iT) / 2)}{\Gamma ((1 - H + iT) / 2)} \notag \\
&<<_{H, n_{K}}  {\left( \frac{{|T|}^{H / 2}}{{|T|}^{- (1 + H) / 2}} \right)}^{a (\chi)}
{\left( \frac{{|T|}^{(H + 1) / 2}}{{|T|}^{- H / 2}} \right)}^{b (\chi)} \notag \\
&<<_{H, n_{K}} {|T|}^{n_{K} (1/2 + H)} \notag \\
&< {(1 + |T|)}^{n_{K} (1 / 2 + H)} \notag
\end{align}
where we apply the Stirling (Lemma ~\ref{T:403}) in the third row
since the poles of the gamma factors are non-positive integers.
For $0 \leq |t| \leq 1$, by continuity, we get $L_{\infty} (s, \chi)$
is bounded. Hence we get (1).
\medskip

(2) Now fix $H, \epsilon' > 0$. Let $s = -H + i T$, then have
\[
G_{1} = \prod_{v \in S} \frac{1 - a_{v}^{-1} q_{v}^{-s}}{1 - a_{v} q_{v}^{1-s}}
\]
where $a_{v} = \chi_{v} (\pi_{v})$ while $\pi_{v}$ the uniformizer of $K_{v}$.
Hence
\begin{align}
|G_{1}| &<= \prod_{v \in S} (1 + q_{v}^{H}) |L_{S} (1 - s, \chi)| \notag \\
     &\leq N_{S}^{H} 2^{|S|} \zeta_{K, S} (1 + H) \notag \\
     &<<_{H, n_{K}} N_{S}^{H} 2^{|S|} \notag
\end{align}

\medskip

It suffices to show that $2^{|S|} <<_{H, n_{K}, \epsilon'} N_{S}^{\epsilon'}$.
In fact, let $b = [2^{\epsilon'}] + 1$, then $2^{|S|} < 2^{n_{K} b} N_{S}^{\epsilon'}$
as when $\fP_{v} > b$ then $\fP_{v}^{\epsilon'} > 2$.

\medskip

Done.

\qedsymbol

\medskip

\emph{Proof of Theorem ~\ref{T:401}, Cont.\ }:

\medskip

Now
\begin{align}
S (X, \chi, S) &=
\frac{1}{2 \pi i} \int_{(-H)} X^{s} W (s) W (\chi) {A (\chi)}^{1/2 - s} G_{0} (s) G_{1} (s)
        L^{S} (1 - s, \chi^{-1}) d\,s \notag \\
&=
\frac{1}{2 \pi i} \int_{(-H)} X^{s} W (s) W (\chi) {A (\chi)}^{1/2 - s} G_{0} (s) G_{1} (s)
        \sum_{n = 1}^{\infty} \frac{a_{\chi, S} (n)}{n^{1 + H}}  d\,s \notag \\
&=
\frac{1}{2 \pi i} \sum_{n = 1}^{\infty} \frac{a_{\chi, S} (n)}{n^{1 + H}}
        \int_{(-H)} X^{s} W (s) W (\chi) {A (\chi)}^{1/2 - s} G_{0} (s) G_{1} (s)
        n^{s + H} d\,s \notag
\end{align}

Here the interchange of the sum and the integral is guaranteed by the absolute convergence
of the Dirichlet series, rapid decay of $W (s)$. Applying Lemma ~\ref{T:402},
we have
\begin{align}
S &(X, \chi, S) <<_{H, K} L^{S} (1 + H, \chi)
\int_{(-H)} \left| X^{s} W (s) W (\chi) {A (\chi)}^{1/2 - s} G_{0} (s) G_{1} (s)
         d\,s \right| \notag \\
&<<_{H, K} {\zeta (1 + H)}^{n_{K}}
X^{-H} {N (\chi)}^{1/2 + H} N_{S}^{H} \int_{(-H)}
    W (s) {(1 + |t|)}^{(1/2 + H) n_{K}} d\,s \notag \\
&<<_{H, K} X^{-H} {N (\chi)}^{1/2 + H} N_{S}^{H} \notag
\end{align}

\medskip

To establish the theorem, we need to bound $S (X, \chi, S)$
below. Now assume that $\chi_{v} = 1$ for all $v \notin S$ such that $N (\fP_{v}) \leq 3 X$.
Then
\[
S (X, \chi, S) = S (X, 1, S)
\]
Note that
$S (X, 1, S)$ is greater than $e^{-1}$ multiples of the primes $\fP_{v}$ of $K$
of degree $1$ outside $S$ such that $X \leq \fP_{v} \leq 2 X$,
By the prime number theorem, for some $A > 0$ depending on $K$, when $X > A$,
$S (X, 1, S) > X / (4 \log (X)) - |S|$. Hence
when $X > 4 A {|S|}^{2}$, $S (X, 1, S) > X / (8 \log (X))$.

\medskip

Then we have, when $X > 4 A {|S|}^{2}$,
\[
X / (4 \log (X)) < S (X, 1, S) = S (X, \chi, S) < C' X^{-H} {N (\chi)}^{1/2 + H} N_{S}^{H}
\]
for some $C'$ depending on $H$ and $K$.

\medskip

Therefore,
\begin{align}
X^{H + 1 - \epsilon} &<<_{\epsilon} X^{H + 1} {\log (X)}^{-1} \notag \\
&<<_{H, K} {N (\chi)}^{1/2 + H} N_{S}^{H} \notag
\end{align}
and thus
\[
X <<_{\epsilon, H, K} {N (\chi)}^{\frac{H + 1/2}{H + 1 - \epsilon}} N_{S}^{\frac{H}{H + 1 - \epsilon}}
\]

\medskip

Hence when $H$ is sufficiently small,
\[
X < C {N (\chi)}^{\epsilon + 1/2} N_{S}^{\epsilon}
\]
where $C$ depends on $\epsilon$ and $K$.
Thus, we have
\begin{align}
X &< {\rm Max} (C {N (\chi)}^{\epsilon + 1/2} N_{S}^{\epsilon}, 4 A {|S|}^{2}) \notag \\
&<<_{\epsilon, K} {N (\chi)}^{\epsilon + 1/2} N_{S}^{\epsilon} \notag
\end{align}

\medskip

Hence Theorem ~\ref{T:401} follows.

\medskip

\qedsymbol

\bigskip

\subsection{$S$-version of Multiplicity One for $GL (n)$}

\medskip \hspace{2ex} \medskip

In fact, such method can be applied to prove the following theorem.
Here the analytic conductor is defined as in \cite{H-R95}, \cite{Brumley06}, \cite{LW2009}.
See \cite{Wang2013-2}.

\medskip

\begin{theorem} \label{T:404}
Let $\pi = \otimes_{v} \pi_{v}$ and $\pi' = \otimes_{v} \pi'_{v}$ be two cuspidal
automorphic representations of $\GL_{d} (\bA_{K})$ ($d \geq 2$), with analytic conductors $\leq Q$,
and let $S$ be a finite set of places of $K$.

\medskip

Then if
Given an arbitrary $H > 0$, then there exists
a constant $C = C (H, d, K, \epsilon)$ such that, if $\pi_{v} \cong \pi'_{v}$
for all places $v$ outside $S$ with the norm
\[
N (\fP_{v}) \leq C {Q}^{2 d + \epsilon} N_{S}^{d^{2} + \epsilon}
\]
Then $\pi \cong \pi'$.
\end{theorem}

\medskip

\qedsymbol

\end{document}